\documentclass[reqno,10pt]{article}
\usepackage{graphicx}
\usepackage{amsmath}
\usepackage{amssymb}
\usepackage{amsthm}
\usepackage{enumitem}
\usepackage{bbm}
\usepackage[l2tabu,orthodox]{nag}
\usepackage[all,error]{onlyamsmath}
\usepackage[strict]{csquotes}
\usepackage{hyperref}
\usepackage{float}
\usepackage{advdate}

\newtheorem{theorem}{Theorem}
\newtheorem{lemma}[theorem]{Lemma}
\newtheorem{corollary}[theorem]{Corollary}
\theoremstyle{definition}
\newtheorem{definition}[theorem]{Definition}
\newtheorem{remark}[theorem]{Remark}

\numberwithin{equation}{section}
\numberwithin{theorem}{section}

\newenvironment{OMabstract}{\noindent\textbf{Abstract.} }{\medskip}
\newenvironment{OMsubjclass}{\noindent\textbf{Mathematics Subject Classification (2020):} }{\medskip}
\newenvironment{OMkeywords}{\noindent\textbf{Keywords:}  }{\medskip}

\DeclareMathOperator{\dist}{dist}
\DeclareMathOperator{\sinc}{sinc}

\begin{document}

\author{Mariusz Niwi\'nski}
\title{Spectral Theory and Continuum Norm Resolvent Convergence of a One-Sided Discrete Dirac-Type Operator}
 \date{October 28, 2026}
\maketitle

\begin{OMabstract}

	We investigate a one-dimensional discrete Dirac-type operator $D_d$ acting on the Hilbert space $\big(\ell^2(\mathbb{Z};\mathbb{C})\big)^2$, obtained from the classical continuous free Dirac operator by replacing the spatial derivative with the one-sided backward difference operator $I - S_1$, where $S_1$ denotes the right shift on $\ell^2(\mathbb{Z};\mathbb{C})$. We establish the fundamental operator-theoretic and spectral properties of $D_d$: unlike symmetric or supersymmetric lattice discretizations, $D_d$ is bounded and non-normal for every non-zero rest mass $m > 0$, with self-commutator $[D_d, D_d^\dag] = -2mc^3(-\Delta_{\mathrm{disc}})\sigma_2$, yet it is quadratically normal. We compute the exact operator norm, determine the exact spectral radius in two mass regimes, and show via Fourier--Laurent analysis that the spectrum $\sigma(D_d)$ is purely continuous and lies on a two-component algebraic curve in the complex plane undergoing a topological transition at $mc^2 = 2c$. Finally, introducing a lattice spacing $h > 0$ and embedding $\big(\ell^2(h\mathbb{Z};\mathbb{C})\big)^2$ isometrically into $\big(L^2(\mathbb{R};\mathbb{C})\big)^2$ via Whittaker--Shannon band-limited interpolation, we prove that the scaled non-normal discrete Dirac-type operator $D_{d,h}$ converges to the continuous self-adjoint Dirac operator $D$ in the norm resolvent sense at the optimal rate $\mathcal{O}(h)$ as $h \to 0^+$, without requiring frequency cutoffs.
\end{OMabstract}

\begin{OMkeywords}
	discrete Dirac operator, shift operator, non-normal operator, quadratic normality, norm resolvent convergence, continuum limit.
\end{OMkeywords}

\begin{OMsubjclass}
	47B39, 47B37, 47A10, 81Q10.
\end{OMsubjclass}


\section{Introduction and preliminaries}

A classical one-dimensional Dirac operator for free electrons is the differential operator $D\colon\big(L^2(\mathbb{R};\mathbb{C})\big)^2 \supset \big(W^{1,2}(\mathbb{R};\mathbb{C})\big)^2 \to \big(L^2(\mathbb{R};\mathbb{C})\big)^2$, given by
\begin{equation}\label{eq:cont_dirac}
	D = -ic\begin{bmatrix}
		0 & 1 \\
		1 & 0
	\end{bmatrix}\frac{d}{d x} + mc^2\begin{bmatrix}
		1 & 0 \\
		0 & -1
	\end{bmatrix}I,
\end{equation}
where $c > 0$ denotes the speed of light and $m \ge 0$ denotes the particle rest mass (see \cite{Thaller1992}). Defined on the direct sum of Sobolev spaces $\mathcal{D}(D) = \big(W^{1,2}(\mathbb{R};\mathbb{C})\big)^2$, the operator $D$ is an unbounded, closed, self-adjoint operator with purely absolutely continuous spectrum
\[
	\sigma(D) = (-\infty, -mc^2] \cup [mc^2, \infty).
\]

One can define a natural discrete analogue of the operator \eqref{eq:cont_dirac} on a one-dimensional lattice and obtain a discrete one-dimensional Dirac-type operator $D_d\colon \big(\ell^2(\mathbb{Z};\mathbb{C})\big)^2 \to \big(\ell^2(\mathbb{Z};\mathbb{C})\big)^2$, given by the formula
\begin{equation}\label{eq:disc_dirac}
	D_d = -ic\begin{bmatrix}
		0 & 1 \\
		1 & 0
	\end{bmatrix}(I - S_1) + mc^2\begin{bmatrix}
		1 & 0 \\
		0 & -1
	\end{bmatrix}I,
\end{equation}
where $S_1\colon \ell^2(\mathbb{Z};\mathbb{C}) \to \ell^2(\mathbb{Z};\mathbb{C})$ is the right shift operator, that is, for any sequence $x = (x_n)_{n\in\mathbb{Z}} \in \ell^2(\mathbb{Z};\mathbb{C})$, we have $(S_1x)_n = x_{n-1}$. Clearly, $S_1$ is a unitary operator on $\ell^2(\mathbb{Z};\mathbb{C})$ whose adjoint (and inverse) is the left shift $(S_1^\dag x)_n = x_{n+1}$.

In the mathematical physics literature, lattice discretizations of the Dirac operator \eqref{eq:cont_dirac} have been extensively studied in the context of spectral theory, supersymmetric quantum mechanics, and continuum limits \cite{CorneanGardeJensen2022,deOliveiraPrado2005,Rabin1982}. When replacing the continuum derivative $\frac{d}{dx}$ on a lattice, one encounters a fundamental structural dichotomy:
\begin{enumerate}[label={\textup{(\roman*)}}, widest=iii, leftmargin=*]
	\item If operator $\frac{d}{dx}$ is replaced by the skew-adjoint symmetric central difference $\frac{1}{2}(S_1^\dag - S_1)$, the resulting operator $D_{\mathrm{sym}}$ is self-adjoint, but its Fourier symbol $i\sin k$ vanishes both at the origin $k = 0$ and at the boundary of the Brillouin zone $k = \pm\pi$. This produces a spurious second Dirac cone at $k = \pm\pi$, known as \emph{fermion doubling} \cite{NielsenNinomiya1981,Rabin1982}. Because of the high-frequency doubler mode, $D_{\mathrm{sym}}$ fails to converge to $D$ in the full-space norm resolvent topology as the lattice spacing $h \to 0^+$ unless an auxiliary projection cutoff onto low frequencies is imposed \cite{CorneanGardeJensen2022}.
	\item If one instead places $-ic(I - S_1)$ in the lower-left off-diagonal entry and its adjoint $ic(I - S_1^\dag)$ in the upper-right entry, one obtains a self-adjoint supersymmetric discretization \cite{deOliveiraPrado2005}, which implicitly modifies the Pauli matrix structure of \eqref{eq:cont_dirac} by introducing a Wilson-type mass term proportional to $\sigma_2(2I - S_1 - S_1^\dag)$.
\end{enumerate}

In contrast, the operator $D_d$ defined in \eqref{eq:disc_dirac} preserves the exact Pauli matrix structure of \eqref{eq:cont_dirac} while employing the one-sided difference operator $\nabla_- := I - S_1$ in both off-diagonal components. As a consequence, $D_d$ exhibits a distinctive combination of operator-theoretic features:
\begin{enumerate}[label={\textup{(\arabic*)}}, widest=3, leftmargin=*]
	\item For $m > 0$, $D_d$ is neither self-adjoint nor normal, yet its square $D_d^2$ is a normal scalar-block operator (that is, $D_d$ is quadratically normal).
	\item The Fourier symbols of the one-sided difference $I - S_1$ are $1 - e^{-ik} = 2ie^{-ik/2}\sin(k/2)$, which vanishes on the torus $\mathbb{T} = \mathbb{R}/2\pi\mathbb{Z}$ \emph{only} at $k = 0$. Thus $D_d$ is completely free of fermion doubling at $k = \pm\pi$, at the cost of its spectrum becoming a two-component curve in the complex plane $\mathbb{C}$.
	\item When scaled by a lattice spacing $h > 0$ and embedded isometrically into the continuum Hilbert space $\big(L^2(\mathbb{R};\mathbb{C})\big)^2$ via Whittaker--Shannon interpolation, the absence of fermion doubling provides global quadratic coercivity across the entire Brillouin zone $[-\pi/h, \pi/h]$. This allows us to prove unconditional $\mathcal{O}(h)$ norm resolvent convergence of $D_{d,h}$ to the continuous self-adjoint Dirac operator $D$ on the full Hilbert space.
\end{enumerate}

The paper is organized as follows. In Section \ref{sec:basic}, we analyze the basic operator-theoretic properties of $D_d$, including its Cartesian decomposition, self-commutator, quadratic normality, exact operator norm, and exact spectral radius. In Section \ref{sec:spectrum}, we determine the exact complex spectrum of $D_d$ and classify its geometric regimes. In Section \ref{sec:continuum}, we construct the Whittaker--Shannon isometric embedding and prove $\mathcal{O}(h)$ norm resolvent convergence in the continuum limit $h \to 0^+$.

\medskip
\noindent\textbf{Notation.} Throughout the paper, we denote by $\mathcal{H} := \big(\ell^2(\mathbb{Z};\mathbb{C})\big)^2 \cong \ell^2(\mathbb{Z};\mathbb{C}^2)$ the discrete Hilbert space equipped with the standard inner product
\[
\begin{split}
	&\langle u, v \rangle_{\mathcal{H}} = \langle u_1, v_1 \rangle_{\ell^2} + \langle u_2, v_2 \rangle_{\ell^2} = \sum_{n\in\mathbb{Z}} \big(u_{1,n}\overline{v_{1,n}} + u_{2,n}\overline{v_{2,n}}\big), \\
    &\qquad u = (u_1, u_2)^T, \; v = (v_1, v_2)^T \in \mathcal{H},
\end{split}
\]
and by $\boldsymbol{\mathcal{B}}(\mathcal{H})$ the $C^*-$algebra of bounded linear operators on $\mathcal{H}$. We write the standard Pauli matrices as
\[
	\sigma_1 = \begin{bmatrix} 0 & 1 \\ 1 & 0 \end{bmatrix}, \qquad
	\sigma_2 = \begin{bmatrix} 0 & -i \\ i & 0 \end{bmatrix}, \qquad
	\sigma_3 = \begin{bmatrix} 1 & 0 \\ 0 & -1 \end{bmatrix}, \qquad
	\sigma_0 = \begin{bmatrix} 1 & 0 \\ 0 & 1 \end{bmatrix}.
\]
In block-operator form on $\mathcal{H} = \ell^2(\mathbb{Z};\mathbb{C}) \oplus \ell^2(\mathbb{Z};\mathbb{C})$, the operator $D_d$ in \eqref{eq:disc_dirac} acts as
\begin{equation}\label{eq:block_Dd}
	D_d = -ic\sigma_1(I - S_1) + mc^2\sigma_3 I = \begin{bmatrix}
		mc^2 I & -ic(I - S_1) \\
		-ic(I - S_1) & -mc^2 I
	\end{bmatrix}.
\end{equation}

\section{Basic operator-theoretic properties of \texorpdfstring{$D_d$}{Dd}}\label{sec:basic}

We begin by computing the adjoint and the Cartesian decomposition of $D_d$ into its real and imaginary self-adjoint parts. Recall that the standard one-dimensional discrete Laplacian on $\ell^2(\mathbb{Z};\mathbb{C})$ is the negative semidefinite self-adjoint operator $\Delta_{\mathrm{disc}} := S_1 + S_1^\dag - 2I$, so that
\begin{equation}\label{eq:laplacian}
	-\Delta_{\mathrm{disc}} = 2I - S_1 - S_1^\dag = I - S_1 - S_1^\dag + S_1S_1^
    \dag = (I - S_1)(I - S_1^\dag).
\end{equation}

\begin{lemma}\label{lem:cartesian}
	Let $c > 0$ and $m \ge 0$. The operator $D_d$ defined in \eqref{eq:disc_dirac} is a bounded linear operator on $\mathcal{H}$ whose adjoint is given by
	\begin{equation}\label{eq:adjoint}
		D_d^\dag = ic\sigma_1(I - S_1^\dag) + mc^2\sigma_3 I = \begin{bmatrix}
			mc^2 I & ic(I - S_1^\dag) \\
			ic(I - S_1^\dag) & -mc^2 I
		\end{bmatrix}.
	\end{equation}
	Moreover, operator $D_d$ admits the Cartesian decomposition $D_d = \operatorname{Re}(D_d) + i\operatorname{Im}(D_d)$, where the self-adjoint components $\operatorname{Re}(D_d) = \frac{1}{2}(D_d + D_d^\dag)$ and $\operatorname{Im}(D_d) = \frac{1}{2i}(D_d - D_d^\dag)$ take the form
	\begin{align}
		\operatorname{Re}(D_d) &= -ic\sigma_1\left(\frac{S_1^\dag - S_1}{2}\right) + mc^2\sigma_3 I, \label{eq:Re_Dd}\\
		\operatorname{Im}(D_d) &= -\frac{c}{2}\sigma_1(2I - S_1 - S_1^\dag) = -\frac{c}{2}\sigma_1(-\Delta_{\mathrm{disc}}). \label{eq:Im_Dd}
	\end{align}
	In particular, since $c > 0$, $D_d$ is never self-adjoint.
\end{lemma}

\begin{proof}
	Since operator $S_1$ is unitary on the space $\ell^2(\mathbb{Z};\mathbb{C})$, it clearly follows that $\|I - S_1\|_{\boldsymbol{\mathcal{B}}(\ell^2(\mathbb{Z};\mathbb{C}))} \le \|I\|_{\boldsymbol{\mathcal{B}}(\ell^2(\mathbb{Z};\mathbb{C}))} + \|S_1\|_{\boldsymbol{\mathcal{B}}(\ell^2(\mathbb{Z};\mathbb{C}))} = 2$. So $I - S_1 \in \boldsymbol{\mathcal{B}}\big(\ell^2(\mathbb{Z};\mathbb{C})\big)$, and hence $D_d \in \boldsymbol{\mathcal{B}}(\mathcal{H})$. Using the self-adjointness of the Pauli matrices $\sigma_1, \sigma_3$ and $(-ic)^\dag = ic$, we obtain
	\[
		D_d^\dag = \big(-ic\sigma_1(I - S_1) + mc^2\sigma_3 I\big)^\dag = ic\sigma_1(I - S_1^\dag) + mc^2\sigma_3 I,
	\]
	which coincides with \eqref{eq:adjoint}. Adding and subtracting $D_d$ and $D_d^\dag$ yields
	\[
    \begin{split}
		\operatorname{Re}(D_d) &= \frac{D_d + D_d^\dag}{2} = -\frac{ic}{2}\sigma_1\big((I - S_1) - (I - S_1^\dag)\big) + mc^2\sigma_3 I \\
        &= -ic\sigma_1\left(\frac{S_1^\dag - S_1}{2}\right) + mc^2\sigma_3 I
    \end{split}
	\]
	and
	\[
		\operatorname{Im}(D_d) = \frac{D_d - D_d^\dag}{2i} = -\frac{c}{2}\sigma_1\big((I - S_1) + (I - S_1^\dag)\big) = -\frac{c}{2}\sigma_1(2I - S_1 - S_1^\dag).
	\]
	Since $2I - S_1 - S_1^\dag = -\Delta_{\mathrm{disc}} \neq 0$ and $c > 0$, we have $\operatorname{Im}(D_d) \neq 0$, hence $D_d \neq D_d^\dag$.
\end{proof}

\begin{remark}
	Identity \eqref{eq:Re_Dd} shows that the Hermitian part $\operatorname{Re}(D_d)$ is precisely the naive symmetric central-difference discrete Dirac operator $D_{\mathrm{sym}}$, while \eqref{eq:Im_Dd} reveals that the skew-Hermitian part $i\operatorname{Im}(D_d)$ acts as a dissipative term proportional to the discrete Laplacian $-\Delta_{\mathrm{disc}}$ coupled to $\sigma_1$.
\end{remark}

We now examine the normality of $D_d$ and its square $D_d^2$.

\begin{theorem}\label{thm:nonnormal}
	For any $c > 0$ and $m \ge 0$, the self-commutator $[D_d, D_d^\dag] := D_d D_d^\dag - D_d^\dag D_d$ is given by
	\begin{equation}\label{eq:self_commutator}
		[D_d, D_d^\dag] = -2mc^3(-\Delta_{\mathrm{disc}})\sigma_2 = -2mc^3(2I - S_1 - S_1^\dag)\begin{bmatrix}
			0 & -i \\
			i & 0
		\end{bmatrix},
	\end{equation}
	and its operator norm is $\big\|[D_d, D_d^\dag]\big\|_{\boldsymbol{\mathcal{B}}(\mathcal{H})} = 8mc^3$. Consequently, $D_d$ is a normal operator on $\mathcal{H}$ if and only if $m = 0$.
\end{theorem}

\begin{proof}
	Using the Cartesian decomposition $D_d = \operatorname{Re}(D_d) + i\operatorname{Im}(D_d)$ from Lemma \ref{lem:cartesian}, we have $[D_d, D_d^\dag] = -2i[\operatorname{Re}(D_d), \operatorname{Im}(D_d)]$. Because $S_1$ and $S_1^\dag$ commute on $\ell^2(\mathbb{Z};\mathbb{C})$, the difference operators $\frac{1}{2}(S_1^\dag - S_1)$ and $-\Delta_{\mathrm{disc}} = 2I - S_1 - S_1^\dag$ commute with each other and with the identity $I$. Since $[\sigma_1, \sigma_1] = 0$ and $[\sigma_3, \sigma_1] = 2i\sigma_2$, we compute
	\begin{equation*}
		\begin{split}
			[\operatorname{Re}(D_d), \operatorname{Im}(D_d)] &= \left[-ic\sigma_1\left(\frac{S_1^\dag - S_1}{2}\right) + mc^2\sigma_3 I, \; -\frac{c}{2}\sigma_1(-\Delta_{\mathrm{disc}})\right]\\
			&= -\frac{mc^3}{2}[\sigma_3, \sigma_1](-\Delta_{\mathrm{disc}}) = -imc^3(-\Delta_{\mathrm{disc}})\sigma_2.
		\end{split}
	\end{equation*}
	Multiplying by $-2i$ establishes equation \eqref{eq:self_commutator}. Since $\sigma_2$ is a unitary operator on $\mathbb{C}^2$ and we have $\|-\Delta_{\mathrm{disc}}\|_{\boldsymbol{\mathcal{B}}(\ell^2(\mathbb{Z};\mathbb{C}))} = 4$, (see \cite{2008Kirsch}, section 3.1) it follows that the norm $\big\|[D_d, D_d^\dag]\big\|_{\boldsymbol{\mathcal{B}}(\mathcal{H})} = 2mc^3\|-\Delta_{\mathrm{disc}}\|_{\boldsymbol{\mathcal{B}}(\ell^2(\mathbb{Z};\mathbb{C}))} = 8mc^3$, which vanishes if and only if $m = 0$.
\end{proof}

Despite the failure of normality for $m > 0$, the anticommutation of the Pauli matrices $\sigma_1$ and $\sigma_3$ endows $D_d$ with quadratic normality.

\begin{theorem}\label{thm:square_normal}
	Let $c > 0$ and $m \ge 0$. The operator $D_d$ satisfies the following structural identities:
	\begin{enumerate}[label={\textup{(\roman*)}}, widest=iii, leftmargin=*]
		\item \textup{(Quadratic normality)} The square of $D_d$ is the block-diagonal operator
		\begin{equation}\label{eq:Dd_squared}
        \begin{split}
			&D_d^2 = \big(m^2c^4 I - c^2(I - S_1)^2\big)\sigma_0 \\
            &= \begin{bmatrix}
				m^2c^4 I - c^2(I - S_1)^2 & 0 \\
				0 & m^2c^4 I - c^2(I - S_1)^2
			\end{bmatrix}.
        \end{split}
		\end{equation}
		In particular, $[D_d^2, (D_d^2)^\dag] = 0$, so $D_d^2$ is a normal operator on $\mathcal{H}$ for all $m \ge 0$ and $c > 0$.
		\item \textup{(Pseudo-Hermiticity / $\operatorname{PT}-$symmetry)} Let $P\colon \ell^2(\mathbb{Z};\mathbb{C}) \to \ell^2(\mathbb{Z};\mathbb{C})$ denote the spatial parity involution $(Px)_n = x_{-n}$. Then $P = P^\dag = P^{-1}$ satisfies $PS_1P = S_1^\dag$, and the unitary involution $J := \sigma_3(P \oplus P)$ on $\mathcal{H}$ satisfies
		\begin{equation}\label{eq:pseudo_hermitian}
			J D_d J^{-1} = D_d^\dag.
		\end{equation}
		\item \textup{(Charge conjugation)} Let $K\colon \mathcal{H} \to \mathcal{H}$ denote componentwise complex conjugation, $(Ku)_n = \overline{u_n}$. Then the antiunitary involution $C := \sigma_1 K$ anticommutes with $D_d$:
		\begin{equation}\label{eq:chiral_conj}
			C D_d C^{-1} = -D_d.
		\end{equation}
	\end{enumerate}
\end{theorem}

\begin{proof}
	\textup{(i)} Since $\sigma_1^2 = \sigma_3^2 = \sigma_0$, $\sigma_1\sigma_3 + \sigma_3\sigma_1 = 0$, and $I - S_1$ commutes with $mc^2 I$, squaring \eqref{eq:block_Dd} yields
	\[
    \begin{split}
		D_d^2 &= (-ic)^2\sigma_1^2(I - S_1)^2 - imc^3(\sigma_1\sigma_3 + \sigma_3\sigma_1)(I - S_1) + m^2c^4\sigma_3^2 I \\ &= \big(m^2c^4 I - c^2(I - S_1)^2\big)\sigma_0.
    \end{split}
	\]
	Since $S_1$ is unitary, $I - S_1$ is normal, and any polynomial in a normal operator is normal; hence $D_d^2$ is normal.

	\textup{(ii)} For any $x \in \ell^2(\mathbb{Z};\mathbb{C})$ and $n \in \mathbb{Z}$, we have $(PS_1Px)_n = (S_1Px)_{-n} = (Px)_{-n-1} = x_{n+1} = (S_1^\dag x)_n$. Since $\sigma_3\sigma_1\sigma_3 = -\sigma_1$ and $\sigma_3^3 = \sigma_3$, conjugating $D_d$ by $J = \sigma_3(P\oplus P)$ gives
	\[
		J D_d J^{-1} = -ic(\sigma_3\sigma_1\sigma_3)(I - PS_1P) + mc^2\sigma_3^3 I = ic\sigma_1(I - S_1^\dag) + mc^2\sigma_3 I = D_d^\dag.
	\]

	\textup{(iii)} Since $S_1$ and $mc^2 I$ have real matrix entries in the canonical basis of $\ell^2(\mathbb{Z};\mathbb{C})$, they commute with $K$. Moreover, $K(-ic\sigma_1)K^{-1} = ic\overline{\sigma_1} = ic\sigma_1$. Because $\sigma_1^3 = \sigma_1$ and $\sigma_1\sigma_3\sigma_1 = -\sigma_3$, we obtain $CD_dC^{-1} = -D_d$.
\end{proof}

We conclude this section by computing the exact operator norm $\|D_d\|_{\boldsymbol{\mathcal{B}}(\mathcal{H})}$ and the exact spectral radius $r(D_d)$, demonstrating a strict norm--spectral-radius gap for all $m, c > 0$.

\begin{theorem}\label{thm:norm_radius}
	For any $c > 0$ and $m \ge 0$, the operator norm of $D_d$ on $\mathcal{H}$ is given by
	\begin{equation}\label{eq:exact_norm}
		\|D_d\|_{\boldsymbol{\mathcal{B}}(\mathcal{H})} = mc^2 + 2c,
	\end{equation}
	whereas the spectral radius $r(D_d) := \sup\{|\lambda| : \lambda \in \sigma(D_d)\}$ is given explicitly by
	\begin{equation}\label{eq:exact_radius}
		r(D_d) = \begin{cases}
			\sqrt{4c^2 - m^2c^4} &\text{ for } 0 \le mc^2 \le \dfrac{2}{\sqrt{3}}c,\\[8pt]
			\dfrac{(mc^2)^{3/2}}{(m^2c^4 - c^2)^{1/4}} &\text{ for } mc^2 > \dfrac{2}{\sqrt{3}}c.
		\end{cases}
	\end{equation}
	In particular, for all $m > 0$ and $c > 0$, we have the strict inequalities
	\begin{equation}\label{eq:norm_gap}
		r(D_d) < \sqrt{m^2c^4 + 4c^2} < mc^2 + 2c = \|D_d\|_{\boldsymbol{\mathcal{B}}(\mathcal{H})}.
	\end{equation}
\end{theorem}
\begin{proof}
	Multiplying \eqref{eq:adjoint} and \eqref{eq:block_Dd} and using \eqref{eq:laplacian}, we compute the positive self-adjoint operator $D_d^\dag D_d$:
	\begin{equation}\label{eq:Dd_star_Dd}
		\begin{split}
			D_d^\dag D_d &= \begin{bmatrix}
				mc^2 I & ic(I - S_1^\dag) \\
				ic(I - S_1^\dag) & -mc^2 I
			\end{bmatrix}
			\begin{bmatrix}
				mc^2 I & -ic(I - S_1) \\
				-ic(I - S_1) & -mc^2 I
			\end{bmatrix}\\
			&= \big(m^2c^4 I + c^2(-\Delta_{\mathrm{disc}})\big)\sigma_0 + mc^3(-\Delta_{\mathrm{disc}})\sigma_2.
		\end{split}
	\end{equation}
	Let $\mathcal{F}\colon \ell^2(\mathbb{Z};\mathbb{C}) \to L^2([-\pi,\pi];\mathbb{C})$ be the unitary discrete Fourier transform defined by $(\mathcal{F}x)(k) = \frac{1}{\sqrt{2\pi}}\sum_{n\in\mathbb{Z}} x_n e^{-ikn}$ for $k \in [-\pi, \pi]$. Under $\mathcal{F}$, the shift $S_1$ is unitarily equivalent to multiplication by $e^{-ik}$, and $-\Delta_{\mathrm{disc}} = 2I - S_1 - S_1^\dag$ is unitarily equivalent to multiplication by $2 - 2\cos k = 4\sin^2(k/2)$. Consequently, $D_d^\dag D_d$ is unitarily equivalent to multiplication on $\big(L^2([-\pi,\pi];\mathbb{C})\big)^2$ by the $2\times 2$ Hermitian matrix symbols
	\[
		\widehat{D_d}(k)^\dag \widehat{D_d}(k) = \big(m^2c^4 + 4c^2\sin^2(k/2)\big)\sigma_0 + 4mc^3\sin^2(k/2)\sigma_2.
	\]
	Since $\sigma_2$ has eigenvalues $\pm 1$, and $\sigma_0$ and $\sigma_2$ commute, the eigenvalues of $\widehat{D_d}(k)^\dag \widehat{D_d}(k)$ for each $k \in [-\pi, \pi]$ are
	\[
		\mu_\pm(k) = m^2c^4 + 4c^2\sin^2(k/2) \pm 4mc^3\sin^2(k/2).
	\]
	Because $m \ge 0$, $c > 0$, and $\sin^2(k/2) \in [0, 1]$ attains its maximum value $1$ at $k = \pm\pi$, we obtain
	\[
    \begin{split}
		&\|D_d\|_{\boldsymbol{\mathcal{B}}(\mathcal{H})}^2 = \|D_d^\dag D_d\|_{\boldsymbol{\mathcal{B}}(\mathcal{H})} \\ &= \sup_{k\in[-\pi,\pi]} \mu_+(k) = m^2c^4 + 4mc^3 + 4c^2 = (mc^2 + 2c)^2,
    \end{split}
	\]
	which proves \eqref{eq:exact_norm}.

	To establish \eqref{eq:exact_radius}, we use Gelfand's spectral radius formula together with Theorem \ref{thm:square_normal}\textup{(i)}. Since $D_d^2$ is a normal operator, $\|D_d^2\|_{\boldsymbol{\mathcal{B}}(\mathcal{H})} = r(D_d^2) = r(D_d)^2$. Under the Fourier transform $\mathcal{F}\oplus\mathcal{F}$, $D_d^2$ is unitarily equivalent to multiplication by $\big(m^2c^4 - c^2(1 - e^{-ik})^2\big)\sigma_0$. Writing $1 - e^{-ik} = 2ie^{-ik/2}\sin(k/2)$ and setting $M := mc^2$ and $s := \sin^2(k/2) \in [0, 1]$, we have $\cos k = 1 - 2s$ and therefore
	\begin{equation*}
		\begin{split}
			\big|m^2c^4 - c^2(1 - e^{-ik})^2\big|^2 &= \big|M^2 + 4c^2 s e^{-ik}\big|^2\\
			&= M^4 + 8M^2c^2 s\cos k + 16c^4 s^2\\
			&= M^4 + 8M^2c^2 s + 16c^2(c^2 - M^2)s^2 =: P(s).
		\end{split}
	\end{equation*}
	Thus $r(D_d) = \max_{s\in[0,1]} P(s)^{1/4}$. Differentiating with respect to $s \in [0, 1]$ gives
	\[
		P'(s) = 8c^2\big(M^2 + 4(c^2 - M^2)s\big).
	\]
	If $0 \le M \le \frac{2}{\sqrt{3}}c$, then $M^2 + 4(c^2 - M^2)s \ge M^2(1 - 4s) + 4c^2 s \ge 0$ on $[0, 1]$, so $P$ is non-decreasing on $[0, 1]$ and attains its maximum at $s = 1$, where $P(1) = (4c^2 - M^2)^2$; hence $r(D_d) = \sqrt{4c^2 - M^2}$. If $M > \frac{2}{\sqrt{3}}c$, then $c^2 - M^2 < 0$ and $P$ is strictly concave on $[0, 1]$ with interior maximum at $s_* = \frac{M^2}{4(M^2 - c^2)} \in (0, 1)$. Evaluating $P(s_*)$ yields
	\[
		P(s_*) = M^4 + \frac{M^4 c^2}{M^2 - c^2} = \frac{M^6}{M^2 - c^2},
	\]
	so $r(D_d) = P(s_*)^{1/4} = M^{3/2}(M^2 - c^2)^{-1/4}$, proving \eqref{eq:exact_radius}. Finally, for all $M, c > 0$ and $s \in [0, 1]$, the triangle inequality gives us $|M^2 + 4c^2 s e^{-ik}| \le M^2 + 4c^2 s \le M^2 + 4c^2$ (where equality cannot hold simultaneously with $s = 1$ and $e^{-ik} = 1$), yielding \eqref{eq:norm_gap}.
\end{proof}

\section{Fourier--Laurent symbol and exact complex spectrum}\label{sec:spectrum}

By means of the discrete Fourier transform $\mathcal{U} := \mathcal{F} \oplus \mathcal{F}\colon \mathcal{H} \to \big(L^2([-\pi,\pi];\mathbb{C})\big)^2$, using techniques similar as in \cite{Thaller1992}, the translation-invariant operator $D_d$ is unitarily equivalent to the multiplication operator $(\mathcal{U} D_d \mathcal{U}^{-1}\widehat{u})(k) = \widehat{D_d}(k)\widehat{u}(k)$, where the $2\times 2$ matrix-valued symbol $\widehat{D_d}\colon [-\pi,\pi] \to M_2(\mathbb{C})$ is given by
\begin{equation}\label{eq:symbol}
	\widehat{D_d}(k) = \begin{bmatrix}
		mc^2 & -ic(1 - e^{-ik}) \\
		-ic(1 - e^{-ik}) & -mc^2
	\end{bmatrix}, \qquad k \in [-\pi, \pi].
\end{equation}

\begin{theorem}\label{thm:spectrum}
	Let $c > 0$ and $m > 0$. The point spectrum $\sigma_p(D_d)$ and the residual spectrum $\sigma_r(D_d)$ are both empty. The spectrum of $D_d$ coincides with its continuous spectrum $\sigma_c(D_d)$ and is given by the complex curve
	\begin{equation}\label{eq:spectrum_curve}
		\sigma(D_d) = \sigma_c(D_d) = \left\{\pm\sqrt{m^2c^4 - c^2(1 - e^{-ik})^2} : k \in [-\pi, \pi]\right\} \subset \mathbb{C},
	\end{equation}
	which is symmetric with respect to both the real and imaginary axes ($\sigma(D_d) = -\sigma(D_d) = \overline{\sigma(D_d)}$). Moreover, the global geometry of $\sigma(D_d)$ is governed by three regimes:
	\begin{enumerate}[label={\textup{(\roman*)}}, widest=iii, leftmargin=*]
		\item \textup{(Heavy-mass regime: $mc^2 > 2c$)} The spectrum $\sigma(D_d) = \Gamma_+ \cup \Gamma_-$ consists of two disjoint simple closed curves $\Gamma_+ \subset \{z \in \mathbb{C} : \operatorname{Re} z > 0\}$ and $\Gamma_- = -\Gamma_+ \subset \{z \in \mathbb{C} : \operatorname{Re} z < 0\}$, separated by a real spectral gap of width $2\sqrt{m^2c^4 - 4c^2}$. The right loop $\Gamma_+$ intersects the real axis at $\sqrt{m^2c^4 - 4c^2}$ (for $k = \pm\pi$) and at $mc^2$ (for $k = 0$).

		\item \textup{(Critical regime: $mc^2 = 2c$)} The spectral gap closes at $k = \pm\pi$, and $\sigma(D_d)$ forms a connected figure-eight curve (lemniscate) passing through the origin $\lambda = 0$ and intersecting the real axis at $\pm 2c$ and $0$.

		\item \textup{(Light-mass regime: $0 < mc^2 < 2c$)} As $k$ ranges over $[-\pi, \pi]$, the two branches connect across the imaginary axis at the points $\pm i\sqrt{4c^2 - m^2c^4}$ (attained at $k = \pm\pi$) and intersect the real axis at $\pm mc^2$ (attained at $k = 0$).
        
	\end{enumerate}
\end{theorem}

\begin{proof}
	For any $\lambda \in \mathbb{C}$ and $k \in [-\pi, \pi]$, the characteristic polynomial of $\widehat{D_d}(k)$ is
	\[
		\det\big(\widehat{D_d}(k) - \lambda \sigma_0\big) = \lambda^2 - \big(m^2c^4 - c^2(1 - e^{-ik})^2\big).
	\]
	Since $k \mapsto \widehat{D_d}(k)$ is continuous and $2\pi$-periodic on $[-\pi, \pi]$, the multiplication operator $\mathcal{U}D_d\mathcal{U}^{-1}$ is invertible on $\big(L^2([-\pi,\pi];\mathbb{C})\big)^2$ if and only if $\det(\widehat{D_d}(k) - \lambda \sigma_0) \neq 0$ for all $k \in [-\pi, \pi]$, which proves \eqref{eq:spectrum_curve}.

	To see that $\sigma_p(D_d) = \emptyset$, suppose $(\widehat{D_d}(k) - \lambda \sigma_0)\widehat{u}(k) = 0$ for Lebesgue-almost every $k \in [-\pi, \pi]$. Then $\det(\widehat{D_d}(k) - \lambda \sigma_0) = 0$ on the support of $\widehat{u}$. However, the Laurent polynomial $z \mapsto \lambda^2 - m^2c^4 + c^2(1 - z^{-1})^2$ is non-constant for $c > 0$ and therefore has at most two roots on the unit circle $|z| = 1$. Thus $\det(\widehat{D_d}(k) - \lambda \sigma_0) = 0$ holds for at most two values of $k \in [-\pi, \pi)$ (a set of Lebesgue measure zero), forcing $\widehat{u} = 0$ in $\big(L^2([-\pi,\pi];\mathbb{C})\big)^2$ and $\sigma_p(D_d) = \emptyset$. By Theorem \ref{thm:square_normal}\textup{(ii)}, $D_d^\dag$ is unitarily equivalent to $D_d$ via $J$, so $\sigma_p(D_d^\dag) = \emptyset$ as well. Since $\sigma_r(D_d) \subset \overline{\sigma_p(D_d^\dag)} = \emptyset$, we conclude that $\sigma(D_d) = \sigma_c(D_d)$.

	Finally, writing $\mathcal{N}(k) := m^2c^4 - c^2(1 - e^{-ik})^2$, we observe that $\mathcal{N}(-k) = \overline{\mathcal{N}(k)}$, so $\sigma(D_d)$ is symmetric across both axes. Moreover, we have $\operatorname{Im}\mathcal{N}(k) = -2c^2\sin k\,(1 - \cos k)$ vanishes on $[-\pi, \pi]$ if and only if $k = 0$ (where $\mathcal{N}(0) = m^2c^4$) or $k = \pm\pi$ (where $\mathcal{N}(\pm\pi) = m^2c^4 - 4c^2$). Evaluating $\pm\sqrt{\mathcal{N}(k)}$ at $k = 0$ and $k = \pm\pi$ according to the sign of $m^2c^4 - 4c^2$ immediately yields regimes \textup{(i)}--\textup{(iii)}.
\end{proof}

The qualitative geometric structure of the complex spectrum $\sigma(D_d)$ across its three mass regimes is illustrated in Figures \ref{fig:spectrum_m4}, \ref{fig:spectrum_m2}, and \ref{fig:spectrum_m1} below.

\begin{figure}[H]
    \centering
    \includegraphics[width=0.55\textwidth]{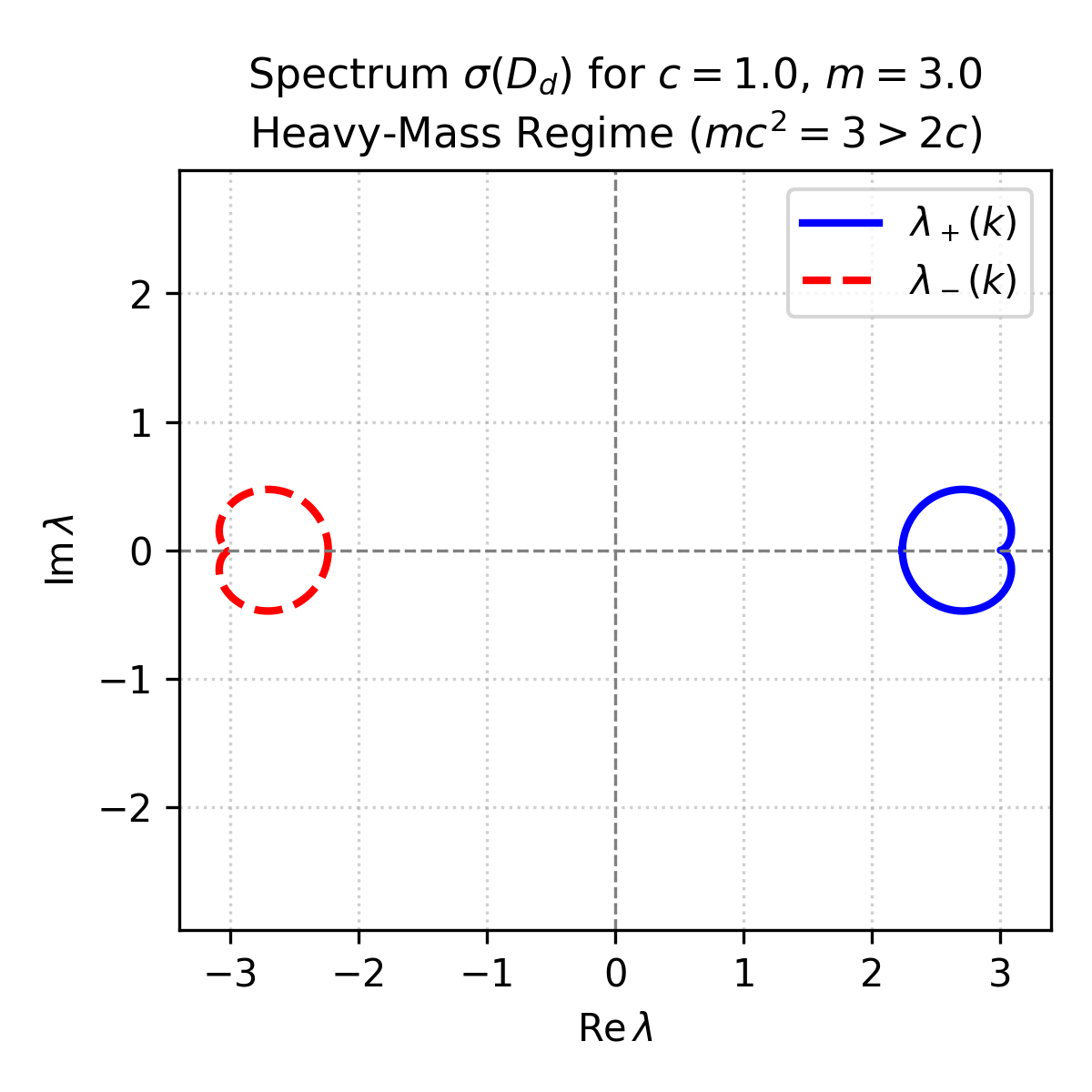}
    \caption{Complex spectrum $\sigma(D_d)$ in the heavy-mass regime ($c=1, m=3$).}
    \label{fig:spectrum_m4}
\end{figure}

\begin{figure}[H]
    \centering
    \includegraphics[width=0.55\textwidth]{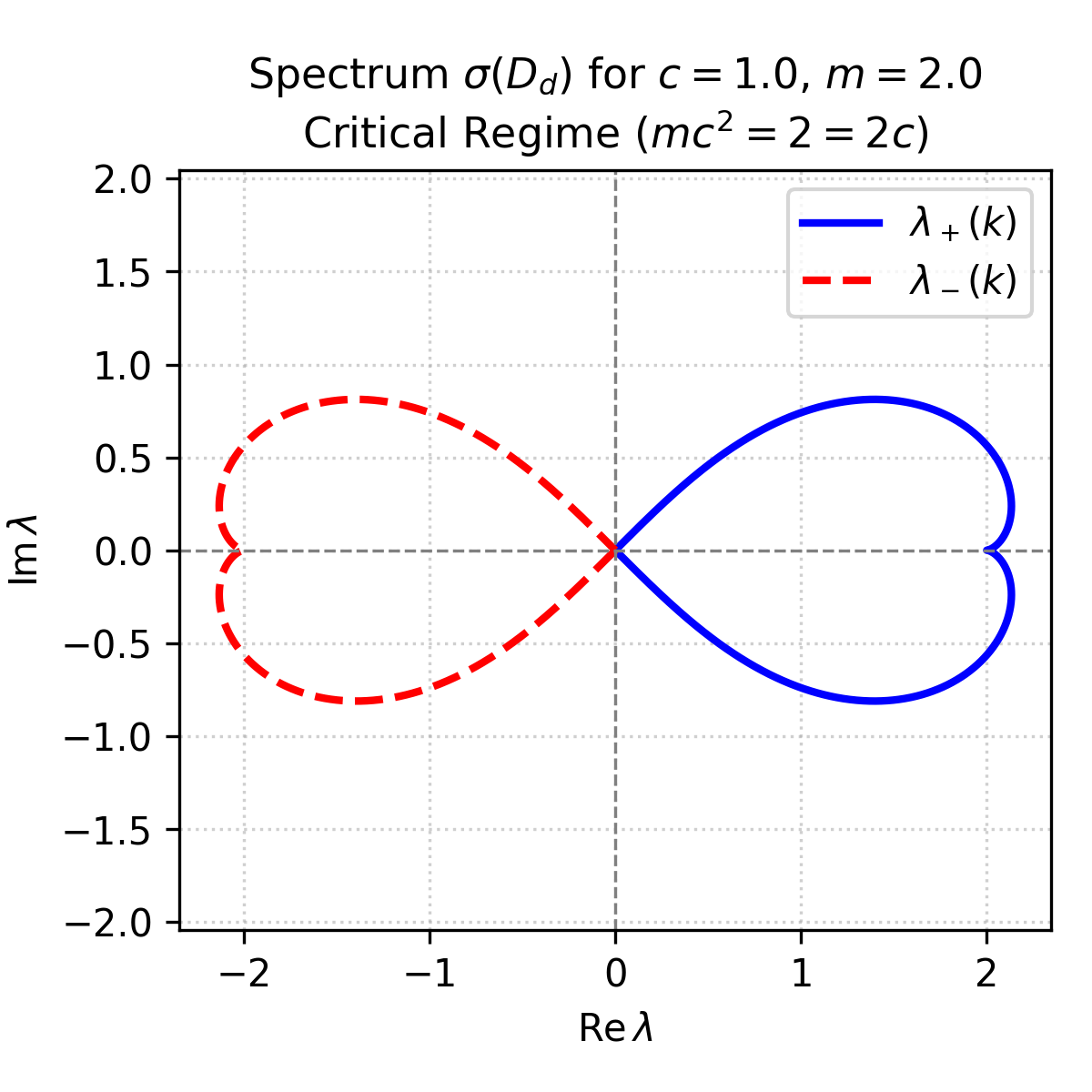}
    \caption{Complex spectrum $\sigma(D_d)$ in the critical regime ($c=1, m=2$).}
    \label{fig:spectrum_m2}
\end{figure}

\begin{figure}[H]
    \centering
    \includegraphics[width=0.55\textwidth]{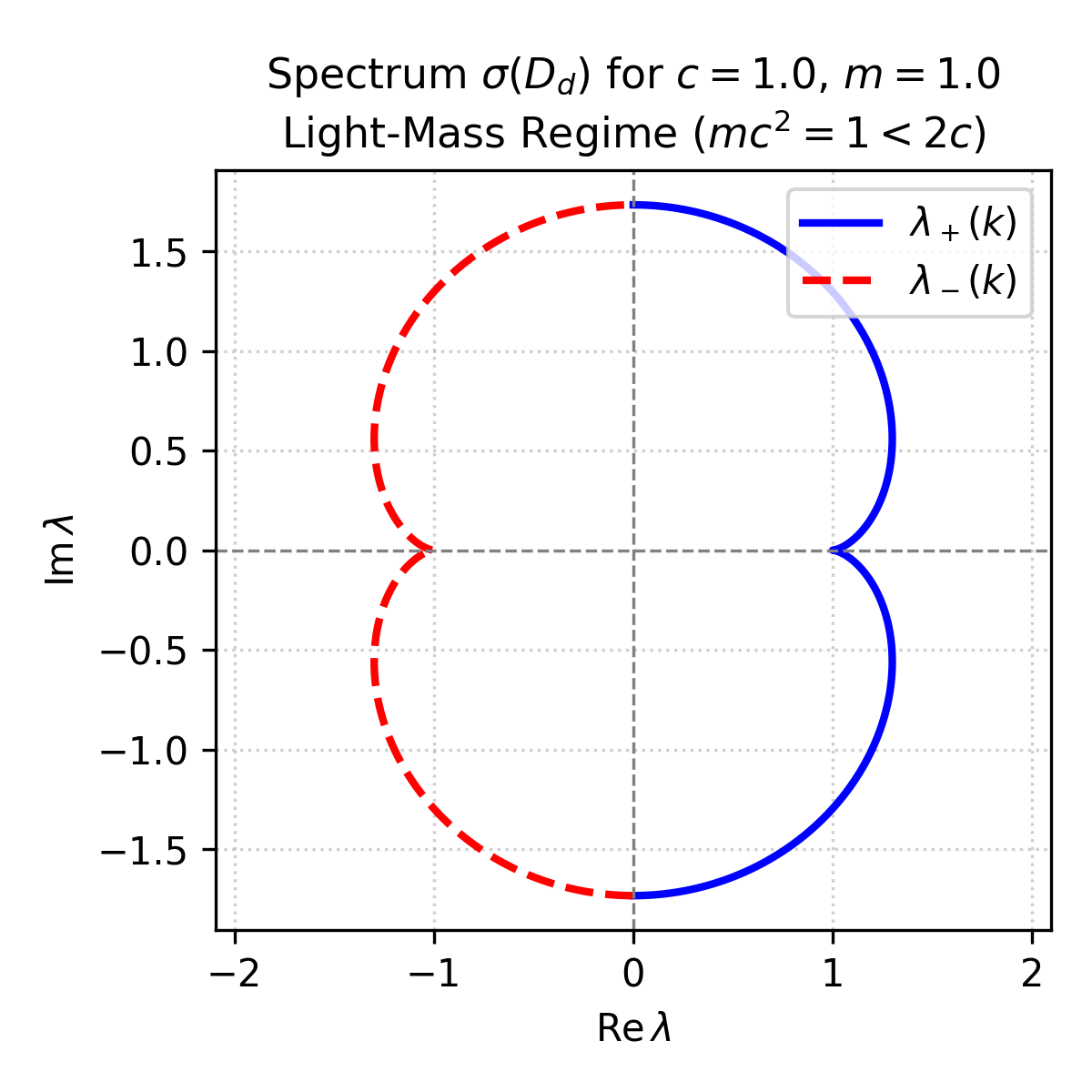}
    \caption{Complex spectrum $\sigma(D_d)$ in the light-mass regime ($c=1, m=1$).}
    \label{fig:spectrum_m1}
\end{figure}

The numerical routines and plotting scripts used to generate these spectral curves are openly available in the GitHub repository \cite{Niwinski2026plots}.

\section{Continuum limit and norm resolvent convergence}\label{sec:continuum}

We now turn to the central question of how the bounded, non-normal discrete operator $D_d$ approximates the unbounded, self-adjoint continuous Dirac operator $D$ in \eqref{eq:cont_dirac} when the lattice spacing tends to zero.

\subsection{Lattice scaling and Whittaker--Shannon isometric embedding}

For each $h > 0$, consider the scaled lattice Hilbert space $\mathcal{H}_h := \big(\ell^2(h\mathbb{Z};\mathbb{C})\big)^2$ endowed with the inner product and norm
\[
	\langle u, v \rangle_{\mathcal{H}_h} := h \sum_{n\in\mathbb{Z}} \big(u_{1,n}\overline{v_{1,n}} + u_{2,n}\overline{v_{2,n}}\big), \qquad \|u\|_{\mathcal{H}_h} = \sqrt{\langle u, u \rangle_{\mathcal{H}_h}},
\]
and define the scaled one-sided discrete Dirac-type operator $D_{d,h}\colon \mathcal{H}_h \to \mathcal{H}_h$ by replacing $I - S_1$ in \eqref{eq:disc_dirac} with the backward difference quotient $\frac{1}{h}(I - S_1)$:
\begin{equation}\label{eq:scaled_Dd}
	D_{d,h} := -\frac{ic}{h}\begin{bmatrix}
		0 & 1 \\
		1 & 0
	\end{bmatrix}(I - S_1) + mc^2\begin{bmatrix}
		1 & 0 \\
		0 & -1
	\end{bmatrix}I = \begin{bmatrix}
		mc^2 I & -\frac{ic}{h}(I - S_1) \\
		-\frac{ic}{h}(I - S_1) & -mc^2 I
	\end{bmatrix}.
\end{equation}

Let $B_h := [-\pi/h, \pi/h]$ denote the scaled Brillouin zone. The scaled discrete Fourier transform $\mathcal{F}_h\colon \ell^2(h\mathbb{Z};\mathbb{C}) \to L^2(B_h;\mathbb{C})$, defined by
\[
	(\mathcal{F}_h f)(p) := \frac{h}{\sqrt{2\pi}}\sum_{n\in\mathbb{Z}} f_n e^{-ipnh}, \qquad p \in B_h,
\]
is unitary with inverse $(\mathcal{F}_h^{-1}\phi)_n = \frac{1}{\sqrt{2\pi}}\int_{-\pi/h}^{\pi/h} \phi(p)e^{ipnh}\,dp$. Likewise, on the Hilbert space $\mathcal{H}_{\mathrm{cont}} := \big(L^2(\mathbb{R};\mathbb{C})\big)^2$, the continuous Fourier--Plancherel transform $\mathcal{F}_{\mathrm{cont}}\colon L^2(\mathbb{R};\mathbb{C}) \to L^2(\mathbb{R};\mathbb{C})$, defined on intersection $L^1(\mathbb{R};\mathbb{C}) \cap L^2(\mathbb{R};\mathbb{C})$ by $(\mathcal{F}_{\mathrm{cont}}\psi)(p) = \frac{1}{\sqrt{2\pi}}\int_{\mathbb{R}} \psi(x)e^{-ipx}\,dx$, is unitary.

\begin{definition}[Whittaker--Shannon isometric embedding]\label{def:embedding}
	For each $h > 0$, let $\iota_h\colon \big(L^2(B_h;\mathbb{C})\big)^2 \hookrightarrow \big(L^2(\mathbb{R};\mathbb{C})\big)^2$ be the isometric zero-extension operator
	\[
		(\iota_h \widehat{u})(p) := \begin{cases}
			\widehat{u}(p) &\text{ for } p \in B_h,\\
			0 &\text{ for } p \in \mathbb{R} \setminus B_h.
		\end{cases}
	\]
	The Whittaker--Shannon interpolation operator $T_h\colon \mathcal{H}_h \to \mathcal{H}_{\mathrm{cont}}$ is defined by
	\begin{equation}\label{eq:Jh_def}
		T_h := (\mathcal{F}_{\mathrm{cont}}\oplus\mathcal{F}_{\mathrm{cont}})^{-1} \iota_h (\mathcal{F}_h \oplus \mathcal{F}_h).
	\end{equation}
\end{definition}

In position space, $T_h$ reconstructs a smooth band-limited $\mathbb{C}^2$-valued function in $\mathcal{D}(D) = \big(W^{1,2}(\mathbb{R};\mathbb{C})\big)^2$ from its lattice values via cardinal $\sinc$ series ($\sinc(t) := \frac{\sin(\pi t)}{\pi t}$ for $t \neq 0$ and $\sinc(0) := 1$):
\[
	(T_h u)(x) = \sum_{n\in\mathbb{Z}} u_n \sinc\!\left(\frac{x - nh}{h}\right), \qquad x \in \mathbb{R}.
\]
Because $\mathcal{F}_h$, $\mathcal{F}_{\mathrm{cont}}$, and $\iota_h$ are isometries, $T_h\colon \mathcal{H}_h \to \mathcal{H}_{\mathrm{cont}}$ is an isometry ($T_h^\dag T_h = I_{\mathcal{H}_h}$), while $P_{B_h} := T_h T_h^\dag$ is the orthogonal projection of $\mathcal{H}_{\mathrm{cont}}$ onto the subspace of $\mathbb{C}^2$-valued functions, band-limited to $B_h$.

\subsection{Global symbol coercivity and uniform resolvent bounds}

Under the Fourier transforms $\mathcal{F}_h \oplus \mathcal{F}_h$ and $\mathcal{F}_{\mathrm{cont}} \oplus \mathcal{F}_{\mathrm{cont}}$, the discrete operator $D_{d,h}$ and the continuous operator $D$ act as multiplication by the $2\times 2$ matrix symbols
\begin{align}
	\widehat{D_{d,h}}(p) &= -ic\xi_h(p)\sigma_1 + mc^2\sigma_3 = \begin{bmatrix}
		mc^2 & -ic\xi_h(p) \\
		-ic\xi_h(p) & -mc^2
	\end{bmatrix}, \qquad p \in B_h, \label{eq:disc_symbol_h}\\
	\widehat{D}(p) &= cp\sigma_1 + mc^2\sigma_3 = \begin{bmatrix}
		mc^2 & cp \\
		cp & -mc^2
	\end{bmatrix}, \qquad p \in \mathbb{R}, \label{eq:cont_symbol}
\end{align}
where $\xi_h\colon B_h \to \mathbb{C}$ is the scaled one-sided difference symbol
\begin{equation}\label{eq:xi_h}
	\xi_h(p) := \frac{1 - e^{-iph}}{h} = \frac{2i}{h}e^{-iph/2}\sin\!\left(\frac{ph}{2}\right).
\end{equation}

The key analytical advantage of the one-sided discretization \eqref{eq:scaled_Dd} over the symmetric central-difference discretization is that $|\xi_h(p)|$ is globally coercive on the entire Brillouin zone $B_h = [-\pi/h, \pi/h]$.

\begin{lemma}\label{lem:xi_estimates}
	For every $h > 0$ and all $p \in B_h = [-\pi/h, \pi/h]$, the symbol $\xi_h(p)$ in \eqref{eq:xi_h} satisfies:
	\begin{align}
		\frac{2}{\pi}|p| \le |\xi_h(p)| &\le |p|, \label{eq:jordan_coercivity}\\
		\big|\xi_h(p) - ip\big| &\le \frac{1}{2}p^2 h, \label{eq:xi_taylor}\\
		\big|\xi_h(p)^2 + p^2\big| &\le |p|^3 h. \label{eq:xi_sq_taylor}
	\end{align}
\end{lemma}

\begin{proof}
	For $p \in [-\pi/h, \pi/h]$, the angle $\theta := ph/2$ belongs to $[-\pi/2, \pi/2]$. Applying Jordan's inequality $\frac{2}{\pi}|\theta| \le |\sin\theta| \le |\theta|$ to $|\xi_h(p)| = \frac{2}{h}|\sin(ph/2)|$ immediately yields \eqref{eq:jordan_coercivity}.

	Next, by Taylor's formula with integral remainder applied to $t \mapsto e^{-ipt}$ on $[0, h]$, we have
	\[
		e^{-iph} = 1 - iph - p^2\int_0^h (h - t)e^{-ipt}\,dt.
	\]
	Dividing by $h$ and taking absolute values gives
	\[
		\big|\xi_h(p) - ip\big| = \frac{p^2}{h}\left|\int_0^h (h - t)e^{-ipt}\,dt\right| \le \frac{p^2}{h}\int_0^h (h - t)\,dt = \frac{1}{2}p^2 h,
	\]
	which proves \eqref{eq:xi_taylor}. Finally, since $|\xi_h(p) + ip| \le |\xi_h(p)| + |p| \le 2|p|$ by \eqref{eq:jordan_coercivity}, multiplying \eqref{eq:xi_taylor} by $|\xi_h(p) + ip|$ gives \eqref{eq:xi_sq_taylor}.
\end{proof}

We next establish that for any fixed spectral parameter $z \in \mathbb{C} \setminus \sigma(D)$, the determinant of $\widehat{D_{d,h}}(p) - z\sigma_0$ satisfies a uniform quadratic lower bound across $B_h$ for all sufficiently small $h > 0$. Let us set
\[
	\mathcal{N}_h(p) := m^2c^4 - c^2\xi_h(p)^2, \qquad E(p) := \sqrt{m^2c^4 + c^2 p^2},
\]
so that $\det(\widehat{D_{d,h}}(p) - z\sigma_0) = z^2 - \mathcal{N}_h(p)$ and $\det(\widehat{D}(p) - z\sigma_0) = z^2 - E(p)^2$.

\begin{lemma}\label{lem:resolvent_bound}
	Let $c > 0$, $m \ge 0$, and $z \in \mathbb{C} \setminus \sigma(D)$, and set $d(z) := \dist(z, \sigma(D)) > 0$. Define the momentum threshold and lattice-spacing threshold
	\begin{equation}\label{eq:thresholds}
		P_0(z) := \frac{\pi}{2c}\sqrt{2\big|m^2c^4 - z^2\big| + \frac{4c^2}{\pi^2}}, \qquad h_0(z) := \frac{d(z)^2}{2c^2 P_0(z)^3}.
	\end{equation}
	Then for every $h \in (0, h_0(z)]$ and all $p \in B_h = [-\pi/h, \pi/h]$, we have
	\begin{equation}\label{eq:det_lower_bound}
		\big|\mathcal{N}_h(p) - z^2\big| \ge \kappa(z)(1 + p^2),
	\end{equation}
	where $\kappa(z) := \min\!\left\{\frac{2c^2}{\pi^2}, \, \frac{d(z)^2}{2(1 + P_0(z)^2)}\right\} > 0$. Consequently, $z \in \rho(D_{d,h})$ for all $h \in (0, h_0(z)]$, and the inverse matrix symbols satisfy
\begin{equation}\label{eq:symbol_res_bounds}
\begin{split}
	& \big\|\big(\widehat{D_{d,h}}(p) - z\sigma_0\big)^{-1}\big\|_{\boldsymbol{\mathcal{B}}(\mathbb{C}^2)} \le \frac{M(z)}{\sqrt{1 + p^2}} \quad (p \in B_h), \\
    & \qquad
	\big\|\big(\widehat{D}(p) - z\sigma_0\big)^{-1}\big\|_{\boldsymbol{\mathcal{B}}(\mathbb{C}^2)} \le \frac{M(z)}{\sqrt{1 + p^2}} \quad (p \in \mathbb{R}),
\end{split}
\end{equation}
	with $M(z) := \frac{1}{\kappa(z)}\sqrt{2(|z| + mc^2)^2 + 2c^2}$.
\end{lemma}

\begin{proof}
	We split $p \in B_h$ into the high-momentum regime $|p| \ge P_0(z)$ and the low-momentum regime $|p| \le P_0(z)$.

	First, suppose $p \in B_h$ satisfies $|p| \ge P_0(z)$. By the lower bound in \eqref{eq:jordan_coercivity} and the definition of $P_0(z)$ in \eqref{eq:thresholds}, we have
	\[
		c^2|\xi_h(p)|^2 \ge \frac{4c^2}{\pi^2}p^2 = \frac{2c^2}{\pi^2}p^2 + \frac{2c^2}{\pi^2}p^2 \ge \frac{2c^2}{\pi^2}(1 + p^2) + \big|m^2c^4 - z^2\big|.
	\]
	Applying the reverse triangle inequality gives
	\[
		\big|\mathcal{N}_h(p) - z^2\big| = \big|c^2\xi_h(p)^2 - (m^2c^4 - z^2)\big| \ge c^2|\xi_h(p)|^2 - \big|m^2c^4 - z^2\big| \ge \frac{2c^2}{\pi^2}(1 + p^2).
	\]

	Second, suppose $|p| \le P_0(z)$. Since $\pm E(p) \in \sigma(D)$, we have $|E(p) - z| \ge d(z)$ and $|E(p) + z| = |-E(p) - z| \ge d(z)$, so $|E(p)^2 - z^2| \ge d(z)^2$. Moreover, by \eqref{eq:xi_sq_taylor}, for all $h \in (0, h_0(z)]$ and $|p| \le P_0(z)$,
	\[
		\big|\mathcal{N}_h(p) - E(p)^2\big| = c^2\big|\xi_h(p)^2 + p^2\big| \le c^2|p|^3 h \le c^2 P_0(z)^3 h_0(z) = \frac{1}{2}d(z)^2.
	\]
	Therefore,
	\[
		\big|\mathcal{N}_h(p) - z^2\big| \ge \big|E(p)^2 - z^2\big| - \big|\mathcal{N}_h(p) - E(p)^2\big| \ge \frac{1}{2}d(z)^2 \ge \frac{d(z)^2}{2(1 + P_0(z)^2)}(1 + p^2),
	\]
	which completes the proof of \eqref{eq:det_lower_bound}.

	To deduce inequalities \eqref{eq:symbol_res_bounds}, observe that $(\widehat{D_{d,h}}(p) - z\sigma_0)(\widehat{D_{d,h}}(p) + z\sigma_0) = (\mathcal{N}_h(p) - z^2)\sigma_0$. Since $\|\sigma_1\|_{\boldsymbol{\mathcal{B}}(\mathbb{C}^2)} = \|\sigma_3\|_{\boldsymbol{\mathcal{B}}(\mathbb{C}^2)} = 1$ and $|\xi_h(p)| \le |p|$, the triangle and Cauchy--Schwarz inequalities give
	\[
		\big\|\widehat{D_{d,h}}(p) + z\sigma_0\big\|_{\boldsymbol{\mathcal{B}}(\mathbb{C}^2)} \le |z| + mc^2 + c|p| \le \sqrt{2(|z| + mc^2)^2 + 2c^2}\sqrt{1 + p^2}.
	\]
	Dividing by $|\mathcal{N}_h(p) - z^2| \ge \kappa(z)(1 + p^2)$ yields the first bound in \eqref{eq:symbol_res_bounds}, and taking $h \to 0^+$ at fixed $p$ yields the second bound.
\end{proof}

\subsection{Optimal \texorpdfstring{$\mathcal{O}(h)$}{O(h)} norm resolvent convergence}

We are now ready to prove our main continuum limit theorem.

\begin{theorem}\label{thm:norm_resolvent}
	Let $c > 0$ and $m \ge 0$, and let $T_h\colon \mathcal{H}_h \to \mathcal{H}_{\mathrm{cont}}$ be the Whittaker--Shannon isometric embedding from Definition \ref{def:embedding}. For every $z \in \mathbb{C} \setminus \sigma(D)$ and all $h \in (0, h_0(z)]$, where $h_0(z)$ is defined in \eqref{eq:thresholds}, we have $z \in \rho(D_{d,h})$ and
	\begin{equation}\label{eq:norm_res_conv}
		\big\|T_h(D_{d,h} - zI)^{-1}T_h^\dag - (D - zI)^{-1}\big\|_{\boldsymbol{\mathcal{B}}(\mathcal{H}_{\mathrm{cont}})} \le C(z)\,h,
	\end{equation}
	where the constant $C(z) > 0$ is given explicitly in terms of $M(z)$ from Lemma \ref{lem:resolvent_bound} by
	\begin{equation}\label{eq:Cz_explicit}
		C(z) := \max\!\left\{\frac{c}{2}M(z)^2, \; \frac{1}{\pi}M(z)\right\}.
	\end{equation}
\end{theorem}

\begin{proof}
	Conjugating the operator difference $T_h(D_{d,h} - zI)^{-1}T_h^\dag - (D - zI)^{-1}$ by the unitary continuous Fourier transform $\mathcal{F}_{\mathrm{cont}}\oplus\mathcal{F}_{\mathrm{cont}}$ on $\mathcal{H}_{\mathrm{cont}}$ and using Definition \ref{def:embedding}, we obtain a multiplication operator on $\big(L^2(\mathbb{R};\mathbb{C})\big)^2$ whose $2\times 2$ matrix multiplier splits orthogonally across $B_h = [-\pi/h, \pi/h]$ and $\mathbb{R} \setminus B_h$:
	\begin{equation}\label{eq:split_sup}
		\begin{split}
			&\big\|T_h(D_{d,h} - zI)^{-1}T_h^\dag - (D - zI)^{-1}\big\|_{\boldsymbol{\mathcal{B}}(\mathcal{H}_{\mathrm{cont}})} = \\ & \max\Big\{\sup_{p\in B_h} \big\|\big(\widehat{D_{d,h}}(p) - z\sigma_0\big)^{-1} - \big(\widehat{D}(p) - z\sigma_0\big)^{-1}\big\|_{\boldsymbol{\mathcal{B}}(\mathbb{C}^2)}, \\&
			\sup_{p\in\mathbb{R}\setminus B_h} \big\|\big(\widehat{D}(p) - z\sigma_0\big)^{-1}\big\|_{\boldsymbol{\mathcal{B}}(\mathbb{C}^2)}\Big\}.
		\end{split}
	\end{equation}

	To bound the first supremum over $p \in B_h$, we apply the second resolvent identity:
	\[
		\big(\widehat{D_{d,h}}(p) - z\sigma_0\big)^{-1} - \big(\widehat{D}(p) - z\sigma_0\big)^{-1} = \big(\widehat{D_{d,h}}(p) - z\sigma_0\big)^{-1}\big(\widehat{D}(p) - \widehat{D_{d,h}}(p)\big)\big(\widehat{D}(p) - z\sigma_0\big)^{-1}.
	\]
	By \eqref{eq:disc_symbol_h}, \eqref{eq:cont_symbol}, and Lemma \ref{lem:xi_estimates}, we have
	\[
		\big\|\widehat{D}(p) - \widehat{D_{d,h}}(p)\big\|_{\boldsymbol{\mathcal{B}}(\mathbb{C}^2)} = c\big|\xi_h(p) - ip\big|\,\|\sigma_1\|_{\boldsymbol{\mathcal{B}}(\mathbb{C}^2)} \le \frac{c}{2}p^2 h.
	\]
	Combining this estimate with the resolvent bounds \eqref{eq:symbol_res_bounds} from Lemma \ref{lem:resolvent_bound} yields, for all $p \in B_h$,
	\begin{equation}\label{eq:interior_est}
		\begin{split}
			&\big\|\big(\widehat{D_{d,h}}(p) - z\sigma_0\big)^{-1} - \big(\widehat{D}(p) - z\sigma_0\big)^{-1}\big\|_{\boldsymbol{\mathcal{B}}(\mathbb{C}^2)} \le \frac{M(z)}{\sqrt{1 + p^2}}\cdot\frac{c}{2}p^2 h\cdot\frac{M(z)}{\sqrt{1 + p^2}}\\
			&= \frac{c}{2}M(z)^2\frac{p^2}{1 + p^2}h \le \frac{c}{2}M(z)^2 h.
		\end{split}
	\end{equation}
	Notice that the quadratic growth $p^2$ of the Taylor remainder is completely absorbed by the product $(1 + p^2)^{-1}$ of the two resolvent decay factors across the entire Brillouin zone $B_h$.

	To bound the second supremum over the high-frequency tail $p \in \mathbb{R} \setminus B_h$, we note that $|p| > \pi/h$ implies $\sqrt{1 + p^2} > |p| > \pi/h$. Applying the second inequality in \eqref{eq:symbol_res_bounds} gives
	\begin{equation}\label{eq:tail_est}
		\sup_{p\in\mathbb{R}\setminus B_h} \big\|\big(\widehat{D}(p) - z\sigma_0\big)^{-1}\big\|_{\boldsymbol{\mathcal{B}}(\mathbb{C}^2)} \le \sup_{|p|>\pi/h} \frac{M(z)}{\sqrt{1 + p^2}} \le \frac{M(z)}{\pi}h.
	\end{equation}
	Substituting \eqref{eq:interior_est} and \eqref{eq:tail_est} into \eqref{eq:split_sup} establishes \eqref{eq:norm_res_conv} with $C(z)$ given by \eqref{eq:Cz_explicit}.
\end{proof}

As a direct consequence, we obtain the precise asymptotic flattening of the complex spectral curve $\sigma(D_{d,h})$ onto the real absolutely continuous spectrum $\sigma(D) = (-\infty, -mc^2] \cup [mc^2, \infty)$ as $h \to 0^+$.

\begin{corollary}\label{cor:spectral_flattening}
	Let $c > 0$ and $m > 0$.
	\begin{enumerate}[label={\textup{(\roman*)}}, widest=iii, leftmargin=*]
		\item For every $h \in (0, \frac{2}{mc})$, $D_{d,h}$ is in the light-mass regime ($mc^2 < 2c/h$). Nevertheless, for each fixed physical momentum $p \in \mathbb{R}$ and $h \to 0^+$, the complex eigenvalues $\lambda_{\pm,h}(p) = \pm\sqrt{\mathcal{N}_h(p)}$ of $\widehat{D_{d,h}}(p)$ satisfy the asymptotic expansion
		\begin{equation}\label{eq:eigenvalue_asymptotics}
			\lambda_{\pm,h}(p) = \pm E(p) \mp \frac{ic^2 p^3}{2E(p)}h + \mathcal{O}(h^2), \qquad E(p) = \sqrt{m^2c^4 + c^2 p^2}.
		\end{equation}
		In particular, the non-Hermitian imaginary part vanishes linearly in $h$: $\operatorname{Im}\lambda_{\pm,h}(p) = \mp\frac{c^2 p^3}{2E(p)}h + \mathcal{O}(h^2)$.
		\item For every compact subset $K \subset \mathbb{C}$, the spectrum $\sigma(D_{d,h})$ converges to $\sigma(D)$ on $K$ in the sense that
		\[
			\sup_{z \in \sigma(D_{d,h}) \cap K} \dist\big(z, \sigma(D)\big) = \mathcal{O}(h) \qquad \text{as } h \to 0^+.
		\]
	\end{enumerate}
\end{corollary}

\begin{proof}
	\textup{(i)} For fixed $p \in \mathbb{R}$ and $h \to 0^+$, the Taylor expansion $\xi_h(p) = ip + \frac{1}{2}p^2 h + \mathcal{O}(h^2)$ implies
	\[
		\mathcal{N}_h(p) = m^2c^4 - c^2\xi_h(p)^2 = E(p)^2 - ic^2 p^3 h + \mathcal{O}(h^2) = E(p)^2\left(1 - \frac{ic^2 p^3}{E(p)^2}h + \mathcal{O}(h^2)\right).
	\]
	Taking the square root via $\sqrt{1 + w} = 1 + \frac{1}{2}w + \mathcal{O}(w^2)$ yields \eqref{eq:eigenvalue_asymptotics}.

	\textup{(ii)} Consider $K \subset \mathbb{C}$ to be a compact set and choose $R > 0$ such that $K \subset \{z \in \mathbb{C} : |z| \le R\}$. If $z \in \sigma(D_{d,h}) \cap K$, then $z = \lambda_{\pm,h}(p)$ for some $p \in B_h$. By the lower bound in \eqref{eq:jordan_coercivity},
	\[
		\frac{4c^2}{\pi^2}p^2 - m^2c^4 \le c^2|\xi_h(p)|^2 - m^2c^4 \le \big|\mathcal{N}_h(p)\big| = |z|^2 \le R^2,
	\]
	so $|p| \le P_R := \frac{\pi}{2c}\sqrt{R^2 + m^2c^4}$. For all $|p| \le P_R$, \eqref{eq:xi_sq_taylor} gives $|\mathcal{N}_h(p) - E(p)^2| \le c^2 P_R^3 h$. Since $\pm E(p) \in \sigma(D)$ and $|\sqrt{\mathcal{N}_h(p)} + E(p)| \ge mc^2 > 0$ for sufficiently small $h > 0$, we obtain
	\[
		\dist\big(z, \sigma(D)\big) \le \big|\lambda_{\pm,h}(p) - (\pm E(p))\big| = \frac{\big|\mathcal{N}_h(p) - E(p)^2\big|}{\big|\sqrt{\mathcal{N}_h(p)} + E(p)\big|} \le \frac{c^2 P_R^3}{mc^2}h,
	\]
	which completes the proof.
\end{proof}

\noindent\textbf{Acknowledgments}\\
\textit{The author would like to thank Professor Petru Cojuhari, for introducing him to this problem, which inspired this research.}

\bigskip


\noindent Mariusz Niwi\'nski\\
mariusz.niwinski@doctoral.uj.edu.pl\\

\noindent {\small
	\noindent Jagiellonian University\\
	Doctoral School of Exact and Natural Sciences\\
	Faculty of Mathematics and Computer Science\\
	ul. \L{}ojasiewicza 6, 30-348 Krak\'ow, Poland
}\bigskip

\end{document}